\documentclass{amsart}
\usepackage{amscd, amsmath, amsthm, amssymb}

\newtheorem{theorem}{Theorem}[section]
\newtheorem{lemma}[theorem]{Lemma}
\newtheorem{proposition}[theorem]{Proposition}
\newtheorem{conjecture}[theorem]{Conjecture}
\newtheorem{corollary}[theorem]{Corollary}

\theoremstyle{definition}     
\newtheorem{definition}[theorem]{Definition}

\newtheorem{claim}[theorem]{Claim}

\theoremstyle{remark}
\newtheorem{remark}[theorem]{Remark}

\numberwithin{equation}{section}

\newcommand{\Alb}{\text{\rm Alb}}
\newcommand{\id}{\text{\rm id}}
\newcommand{\N}{\text{\rm N}}
\newcommand{\NE}{\text{\rm NE}}
\newcommand{\NS}{\text{\rm NS}}

\begin{document}

\title[conjecture of Tits type for
complex varieties] {Conjecture of Tits type for complex varieties
and Theorem of Lie-Kolchin type for a cone}

\author[JongHae Keum, Keiji Oguiso and De-Qi Zhang]{JongHae Keum, Keiji Oguiso, and De-Qi Zhang}
\address{School of Mathematics, Korea Institute for Advanced Study,
Dongdaemun-gu, Seoul 130-722, Korea} \email{jhkeum@kias.re.kr}
\thanks {The first named author was supported by Korea Research Foundation Grant
(KRF-)}

\address{Department of Economics, Keio University, Hiyoshi, Kohoku-ku, Yokohama, Japan and
School of Mathematics, Korea Institute for Advanced Study,
Dongdaemun-gu, Seoul 130-722, Korea }
\email{oguiso@hc.cc.keio.ac.jp}
\thanks {The second named author was supported by JSPS}

\address{Department of Mathematics, National University of Singapore,
2 Science Drive 2, Singapore 117543, Singapore}
\email{matzdq@math.nus.edu.sg}


\begin{abstract} First, we shall formulate and prove Theorem
of Lie-Kolchin type for a cone and derive some algebro-geometric
consequences. Next, inspired by a recent result of Dinh and Sibony
we pose a conjecture of Tits type for a group of automorphisms of
a complex variety and verify its weaker version. Finally, applying
Theorem of Lie-Kolchin type for a cone, we shall confirm the
conjecture of Tits type for complex tori, hyperk\"ahler manifolds,
surfaces, and minimal threefolds.

\end{abstract}

\maketitle


\setcounter{section}{0}
\section{Introduction}

We often study a group $G$ of automorphisms of a projective
variety $X$ through its action on cohomological spaces such as on
the N\'eron-Severi group $\NS(X)$. Let
$$G^{*} = {\rm Im}(G\longrightarrow {\rm GL}\, (\NS(X))).$$
\par
\vskip 4pt Then Tits' Theorem [Ti] implies that

\begin{itemize}
\item [(i)] either $G^{*}$ contains a subgroup isomorphic to
$\mathbf Z * \mathbf Z$ (highly noncommutative);
\item[(ii)] or $G^{*}$ contains a connected solvable subgroup of finite
index.
\end{itemize}

\par
\vskip 4pt Here, {\it connected} means that its Zariski closure in
${\rm GL}\, (\NS(X)_{\mathbf C})$ is connected.

This property already provides a useful tool in studying groups of
automorphisms of projective varieties (see [Og2], [Zh1] for some
concrete applications). Another important feature of $G^{*}$ is
that it preserves the ample cone ${\rm Amp}\, (X)$ and the nef
cone $\overline{{\rm Amp}}\, (X)$, the closure of the ample cone,
both of which encode much information on geometry of $X$.
\par
\vskip 4pt

Keeping these two features in mind, we shall prove first the
following Theorem of Lie-Kolchin type (see also Theorem 2.1 and
Corollary 2.3),
and derive some direct
consequences for groups of automorphisms of projective varieties
or compact K\"ahler manifolds (Theorems 3.1-3.5):

\begin{theorem} \label{theorem:lieint} {\bf (Theorem of Lie-Kolchin type for a cone)}
Let $V$ be a finite dimensional real vector space, and let $C
\not= \{0\}$ be a strictly convex closed cone of $V$. Let $G$ be a
connected $($in the sense above$)$ solvable subgroup of ${\rm
GL}(V)$ such that $G(C) \subseteq C$. Then $G$ has a common
eigenvector in the cone $C$.
\end{theorem}

This theorem is purely algebraic and {\it a priori} has  nothing to do
with algebraic geometry or complex geometry. In our proof, we use
Birkhoff-Perron-Frobenius' Theorem [Bir] and the proof of
Dinh-Sibony [DS] for the case where $G$ is abelian.
\par
\vskip 4pt

In the same paper, Dinh and Sibony proved the following
interesting result:

\begin{theorem} \label{theorem:dsint} $([$DS$])$
Let $M$ be an $n$-dimensional compact K\"ahler manifold and let
$G$ be an abelian subgroup of ${\rm Aut}\, (M)$ such that each
element of $G \setminus \{\id\}$ is of positive entropy. Then $G$
is a free abelian group of rank at most $n -1$. Furthermore, the
rank estimate is optimal.
\end{theorem}

See for instance [ibid] or (4.1) for the notion of entropy of an
automorphism.
\par
\vskip 4pt

In view of Tits' Theorem and Dinh-Sibony's Theorem, it is natural
to pose the following conjecture:

\begin{conjecture} \label{conjecture:titstypeint} {\bf (Conjecture of Tits type)}
Let $X$ be an $n$-dimensional compact K\"ahler manifold or an $n$-dimensional
complex projective variety with at most rational $\mathbf Q$-factorial singularities.
Let $G$ be a subgroup of ${\rm Aut}\, (X)$. Then,
one of the following two assertions holds:

\begin{itemize}
\item[(1)] $G$ contains a subgroup isomorphic to $\mathbf Z * \mathbf Z$.
\item[(2)]
There is a finite index subgroup $G_1$ of $G$ such that the
subset
$$N(G_1) = \{g \in G_1\, \vert\, g\,\ {\rm is\,\ of\,\ null\,\ entropy}\}\,$$
of $G_1$ is a normal subgroup of $G_1$ and the quotient group
$G_1/N(G_1)$ is a free abelian group of rank at most $n-1$.
\end{itemize}
\end{conjecture}

It turns out that the crucial part of the conjecture is the rank
estimate in the statement (2), as a weaker version of the
conjecture (Theorem 4.3) can be easily verified. We shall also
verify the conjecture for some basic cases, say, for surfaces
(Theorem 4.4), hyperk\"ahler manifolds (Theorem 4.5), complex tori
(Theorem 4.6), and minimal threefolds (Theorem 5.1). Except for
the case of complex tori, Theorem of Lie-Kolchin type (Theorem
1.1) plays a crucial role in our proof.

It would be interesting to verify/deny the conjecture for higher
dimensional Calabi-Yau manifolds and manifolds with Kodaira
dimension $-\infty$.

\par \vskip 1pc

{\it Acknowledgement.} We would like to express our thanks to
Professor Yujiro Kawamata for his encouragement and valuable
suggestions. The first two authors would like to thank Department
of Mathematics, National University of Singapore, and staff
members there for financial support and hospitality. The main part of
this work was done during their stay in Singapore.

\par \vskip 1pc

\section{Theorem of Lie-Kolchin type for a cone}

{\bf (2.1)} We call a group  {\it virtually solvable} if it has a
solvable subgroup of finite index. A solvable subgroup $G$ of
${\rm GL}(n, \mathbf C)$ is called {\it connected} if its Zariski
closure $\overline{G}$ in ${\rm GL}(n, \mathbf C)$ is connected.
Note that given a virtually solvable subgroup $G$ of ${\rm GL}(n,
\mathbf C)$, one can always find a connected solvable, finite
index subgroup $H$ of $G$. Indeed, choose a solvable finite index
subgroup  $G_1$ of $G$ and take the identity component
$\overline{G_1}^{0}$ of the Zariski closure $\overline{G_1}$. Then
the group $H := \overline{G_1}^{0} \cap G_1$ is a connected solvable
finite index subgroup of $G$.

\medskip
{\bf (2.2)} Let $V$ be an $r$-dimensional real vector space.
We regard $V$ as an
$r$-dimensional euclidean space with the natural topology. Let $C$
be a subset of $V$. We call $C$ a {\it strictly convex
closed cone of} $V$ if $C$
is closed in $V$, closed under addition and
multiplication by a non-negative scalar, and contains no
1-dimensional linear space.

\medskip
{\bf (2.3)} In this section, we shall prove the following theorem:

\begin{theorem} \label{theorem:lie} {\bf (Theorem of Lie-Kolchin type for a cone)}
Let $V$ be an $r$-dimensional real vector space, and let $C \not=
\{0\}$ be a strictly convex closed cone of $V$. Let $G$ be a
virtually solvable subgroup of ${\rm GL}(V)$ such that $G(C)
\subseteq C$. Then there are a finite-index subgroup $H$ of $G$
and a real vector $v \in C \setminus \{0\}$ such that $H(\mathbf
R_{\ge 0} v) = \mathbf R_{\ge 0} v$. In other words, a suitable
finite-index subgroup of $G$ has a common eigenvector in the cone
$C$. Moreover, if $G$ is connected and solvable, then one may
choose $H = G$.
\end{theorem}

Before starting the proof, we recall several relevant results
and state one special but useful version.

\begin{remark} \label{remark:history}
$ $
\par
(1) Tits' Theorem [Ti] states that any subgroup $G$ of ${\rm
GL}(V_{\mathbf C})$ is either virtually solvable or contains a
free subgroup isomorphic to $\mathbf Z * \mathbf Z$. In the latter
case, the conclusion in Theorem 2.1 is not true in general even in
the case $V=\mathbf Z^r\otimes \mathbf R$ and $G \subset {\rm
GL}\, (\mathbf Z^r)$. In this sense, our restriction on the group
$G$ is reasonable.

(2) The original version of Lie-Kolchin Theorem (see eg. [Hu])
states that any connected solvable subgroup $G$ of ${\rm
GL}(V_{\mathbf C})$ has a common eigenvector. Note that it is, in
general, far from being real, even in the case $V=\mathbf
Z^r\otimes \mathbf R$ and $G \subset {\rm GL}\, (\mathbf Z^r)$.

(3) Birkhoff-Perron-Frobenius' Theorem [Bir] states that any $g
\in {\rm GL}(V)$ such that $g(C) \subset C$ has an eigenvector
$v_{g}$ in $C$. Note that $g$ acts on the real vector space $W =
\langle C \rangle$. So, one can choose $v_{g}$ so that $g(v_{g}) =
\rho(g \vert W) \ v_{g}$, where $\rho(g \vert W)$ is the spectral
radius of $g \vert W$ (see [ibid]).

(4) Dinh-Sibony ([DS] Proposition 4.1) proved Theorem 2.1 when $G$
is abelian. So, our theorem is not only a refinement of
Lie-Kolchin Theorem but also an optimal generalization of
Birkhoff-Perron-Frobenius' Theorem [Bir] and  Dinh-Sibony's
version [DS] in view of Tits' Theorem.
\end{remark}

For algebro-geometric applications, the following arithmetical version
will be most useful:

\begin{corollary} \label{theorem:liekol}
Let $L$ be a free $\mathbf Z$-module of finite rank $r \ge 1$, and
let $C \not= \{0\}$ be a strictly convex closed cone of
$L_{\mathbf R} := L \otimes \mathbf R$. Let $H$ be a connected
solvable subgroup of ${\rm Aut}\, L = {\rm GL}(L)$ such that $H(C)
\subseteq C$. Then there is a real vector $v \in C \setminus
\{0\}$ such that $H(\mathbf R_{\ge 0} v) = \mathbf R_{\ge 0} v$.
\end{corollary}

This version together with the next lemma or its proof
will be frequently used in the sequel.
In this lemma, the condition ``the action of $H$ is
defined over $\mathbf Z$" is very crucial.

\begin{lemma} \label{lemma:disc}
Let $L$ be a free $\mathbf Z$-module of finite rank $r \ge 1$ and
let $H$ be a connected solvable subgroup of ${\rm GL}(L)$. Then, the
following set
$$N(H) := \{h \in H\, \vert\, \rho(h) = 1\}$$
is a normal subgroup of $H$ and $H/N(H)$ is a free abelian group
of rank at most $r-1$. Here, $\rho(h)$ is the spectral radius of
$h$ on $L_{\mathbf C}$.
\end{lemma}

\begin{proof} By Lie-Kolchin Theorem,
$H$ is a subgroup of the
group $T(r)$ of upper triangle matrices of ${\rm GL}(r, \mathbf
C)$ under the identification ${\rm GL}(L_{\mathbf C}) = {\rm
GL}(r, \mathbf C)$ with respect to a suitable basis of $L_{\mathbf
C}$.

Regarding $H \subset T(n)$, we write the $(i,i)$-th entry of the
matrix $h \in H$ by $\chi_{i}(h)$ ($1 \le i \le r$). Note that
$\{\chi_{i}(h)\}_{i=1}^{r}$ is the set of eigenvalues of $h$
(listed with multiplicity) and
$$\chi_{i}(h_{1}h_{2}) =\chi_{i}(h_{1})\chi_{i}(h_{2}),$$
$$\vert \chi_{1}(h) \chi_{2}(h) \cdots \chi_{r}(h) \vert = 1.$$
The last equality holds because $h \in {GL}(L)$ and $L$ is a free
$\mathbf Z$-module. Thus, we have a group homomorphism defined as follows:
$$\chi : \, H \longrightarrow \mathbf R^{r-1} \ =
\ \{(x_{i})_{i=1}^{r} \ \in \ \mathbf R^{r}\, \vert \, \sum_{i=1}^{r} x_{i} = 0\}
\, \subset \, \mathbf R^{r}, $$
$$\chi(h) := (\log \vert \chi_{1}(h) \vert,\  \log \vert \chi_{2}(h)\vert,\
\cdots , \ \log \vert \chi_{r}(h)\vert )\,\, .$$
By the definition of $N(H)$, we have
$$N(H) = {\rm Ker}\, \chi.$$
Thus $N(H)$ is a normal subgroup and its quotient group $H/N(H)
\simeq \chi(H)$ is a subgroup of $\mathbf R^{r-1}$.

It now suffices to show that $\chi(H)$ is discrete in the additive
group $\mathbf R^{r-1}$ with respect to the usual topology. Or
equivalently, it suffices to show that the identity element
$(0)_{i=1}^{r}$ is an isolated point of $\chi(H)$. Here the fact that
$H$ is defined over $\mathbf Z$ is very crucial.

Let $\delta$ be any positive real number. Consider the set of all
elements $h\in H$ satisfying the inequalities
$$\log \vert \chi_{i}(h) \vert < \delta\,\, ,\,\, i = 1, \cdots , r\,\, .$$
For such an $h$, \ $\chi_{i}(h)$ are all bounded, and hence the
coefficients of the characteristic polynomial $\Phi_{h}(x)$ are
all bounded. Since $\Phi_{h}(x) \in \mathbf Z[x]$, there are only
finitely many characteristic polynomials for all such $h$. Thus
the set of vectors $(\log \vert \chi_{i}(h)\vert )_{i=1}^{r}$ for
all such $h$ is also finite. Hence $(0)_{i=1}^{m}$ is an isolated
point of $\chi(H)$.
\end{proof}

In the rest of this section, we shall prove Theorem 2.1. First, we
observe the following:

\begin{lemma} \label{lemma:commut}
Let $H$ be a connected solvable, possibly discrete, subgroup of
${\rm GL}(r, \mathbf C)$. Let $[H,H]$ denote the commutator
subgroup of $H$. Then the eigenvalues of every element of $[H,H]$
are all equal to $1$.
\end{lemma}

\begin{proof} By Lie-Kolchin Theorem, $H$ is a subgroup of the
group $T(r)$ of upper triangle matrices of ${\rm GL}(r, \mathbf
C)$. From this, it is easy to see that the diagonal entries of any
commutator $g_1g_2g_1^{-1}g_2^{-1}$ are all equal to $1$. So are
the diagonal entries of any product of commutators.
\end{proof}

\medskip
{\it Proof of Theorem 2.1.} By replacing $G$ by a connected
solvable subgroup of finite index, we may assume that $G$ itself
is connected and solvable.

We consider the derived series of $G$:

$$G = G^{(0)} > G^{(1)} > \cdots > G^{(k-1)} > A := G^{(k)} > G^{(k+1)}
= \{\ \id\}\,.$$ Here, the series has length $k \ge 0$, $G^{(l)}$
are normal subgroups of $G$ and $[G^{(l)}, G^{(l)}] = G^{(l+1)}$
for $l = 0, 1, \cdots, k$. Note that $A := G^{(k)}$ is an abelian
group.

We shall prove Theorem 2.1 by induction on the length $k$.

If $k = 0$, then $G = A$, i.e. $G$ itself is abelian, and the
result follows from the result of Dinh-Sibony (Remark 2.2(4)).

Assume that $k\ge 1$. Assuming that Theorem 2.1 is true when the
length $\le k-1$, we shall show that Theorem 2.1 is also true when
the length $= k$.

First note that $A\subseteq [G, G] = G^{(1)}$. By Lemma 2.5, the
eigenvalues on $V_{\mathbf C}$ of every element of $A$ are all
equal to $1$.

For each $a\in A$, we write the associated real eigenspace by
$V(a, 1)$:
$$V(a, 1) := \{\, x \in V\, \ \vert\, \  a(x) =
x\, \} \subset V\,\, .$$

Let us consider the linear subspace $W$ of $V$,
$$W=\cap_{a \in A} \, V(a, 1)\,.$$

\begin{claim} \label{claim:W}
$ $
\begin{itemize}
\item[(1)]
$C_W:= C\cap W\not=\{0\}\,\,.$
\item[(2)] For each $g \in G$, $g(W) = W$. In particular,
$G(C_W)\subseteq C_W$.
\end{itemize}
\end{claim}

\begin{proof} (1) Since Theorem 2.1 is true for the abelian
group $A$, there is a common eigenvector in $C$ for $A$ and hence
(1) follows.

\par
(2) Let $0 \ne w \in W$. Since $A$ is a normal subgroup of $G$, we
can write $g^{-1}ag = a_g$ for some $a_g \in A$. Thus, regardless
of the choice of $w$ we have
$$(g^{-1}ag) \ (w) = a_{g}(w) = w\ .$$
Hence
$$a(g(w)) = g(w).$$
So
$$g(W)\, \subseteq \, \cap_{a \in A} \, V(a, 1) \, = W.$$
Since $g$ is invertible, this implies $g(W) = W$.
\end{proof}

From now on, let us consider the action of $G$ on $W$ and on
$C_{W}$. Set
$$G_{W}^{(l)} :=
{\rm Im}\, (G^{(l)} \longrightarrow {\rm GL}(W))\,\, .$$ Then,
$G_{W} := G_{W}^{(0)}$ is connected solvable, and the new series
$$G_{W} = G_{W}^{(0)} \ge G_{W}^{(1)} \ge \cdots \ge G_{W}^{(k-1)} \ge A_{W} :=
G_{W}^{(k)} \ge G_{W}^{(k+1)} = \{\, \id_W\}$$ gives the derived
series, possibly redundant, of $G_{W}$.

Since $W = \cap_{a \in A} \, V(a, 1)$ and hence $A_W = \{\ \id_W\}$,
{\it the new series is in fact redundant}. Now we can apply induction
hypothesis to the series of length $\le k-1$ and we obtain the
result.
\medskip

This completes the proof of Theorem 2.1.

\section{Some algebro-geometric consequences of Theorem of Lie-Kolchin type
for a cone}

In this section, we shall derive some direct consequences of
Theorem of Lie-Kolchin type for a cone (Theorem 2.1 or Corollary 2.3)
for groups of
automorphisms or birational automorphisms of algebraic varieties and compact K\"ahler
manifolds.

\medskip
{\bf (3.1)} Let $f : X \longrightarrow S$ be a projective
surjective morphism of normal varieties defined over a field $k$
with connected fibers of positive dimension. Assume that $S$ is
quasi-projective. Special important cases are the case where $S =
{\rm Spec}\, k$ (abstract case) and the case where $f :=
\Phi_{\vert mK_{X} \vert} : X \longrightarrow S$ is the relatively
minimal Kodaira-Iitaka fibration of $X$.  We note that the
relative canonical class $K_{X/S}$ is trivial in $\NS(X/S)$ (cf.
(3.3)) in the second case.

\medskip
{\bf (3.2)} By ${\rm Bir}\, (X)$ and by ${\rm Aut}\, (X)$, we
denote the group of birational transformations of $X$ and the
group of biregular transformations of $X$ respectively. Let:
$${\rm Bir}(X/S) := \{\varphi \in {\rm Bir}\, (X)\, \vert\, f = f \circ \varphi\}\, ,$$
$${\rm Bir}^{s}(X/S) := \{\varphi \in {\rm Bir}\, (X/S)\, \vert\, \varphi\,\, \
{\rm is\,\ isomorphic\,\ in\,\ codimension\,\ one}\, \}\, ,$$
$${\rm Aut}(X/S) := \{\varphi \in {\rm Aut}\, (X)\, \vert\, f = f \circ \varphi\}\,\, .$$
%
%
%
%
Note that ${\rm Bir}^{s}\,(X) = {\rm Bir}\,(X)$ when $X$ is a (terminal)
minimal model.

\medskip
{\bf (3.3)} Let $\NS(X/S)$ be the set of numerical equivalence
classes of Cartier divisors on $X$ over $S$. Note that $\NS(X/S)$
is the free part of the "N\'eron-Severi" group. It is well known
that $\NS(X/S)$ is a free $\mathbf Z$-module of finite rank ([Kl]).
We define the relative nef cone $\overline{{\rm Amp}}(X/S)$ to
be the closure of the relative ample cone ${\rm Amp}(X/S)$ in
$\NS(X/S)_{\mathbf R}$ (under the usual topology). Similarly, the
relative pseudo-effective cone $\overline{{\rm Big}}(X/S)$ is
defined to be the closure of the relative big cone ${\rm
Big}(X/S)$ in $\NS(X/S)_{\mathbf R}$. Both spaces form closed
convex cones of $\NS(X/S)$. Moreover, since the dimension of fibers
of $f$ is assumed to be positive, both cones are strict. For more
details on these cones and their roles in birational algebraic
geometry, we refer the readers to [Ka2].

\medskip
{\bf (3.4)} Let $G$ be a subgroup of ${\rm Aut}\,(X/S)$. We call
$G$ {\it numerically virtually solvable}  if its image in ${\rm
GL}\, (\NS(X/S))$ is virtually solvable. By Tits' Theorem (Remark
2.2(1)), $G$ satisfies either one of the following two, which are
not in general mutually exclusive:

\begin{itemize}
\item[(i)] $G$ is numerically virtually solvable, or
\item[(ii)] $G$ contains a subgroup isomorphic to $\mathbf Z * \mathbf Z$\, .
\end{itemize}

\medskip
Applying Corollary 2.3 to the action of $G$ on $(\NS(X/S), \,
\overline{{\rm Amp}}(X/S))$, we have the following:

\begin{theorem} \label{theorem:aut}

Let $G$ be a numerically virtually solvable subgroup of ${\rm
Aut}\,(X/S)$. Then there are a finite-index subgroup $G_1$ of $G$
and a nonzero real vector $v \in \overline{{\rm Amp}}(X/S)$ such
that $G_1(\mathbf R_{\ge 0}v) = \mathbf R_{\ge 0}v$. In particular,
if $G$ is a numerically virtually solvable subgroup of ${\rm
Aut}\,(X)$, then a suitable finite-index subgroup of $G$ admits a
common real nef eigenvector.
\end{theorem}

We know that the relative canonical divisor class $K_{X/S}$ is
always fixed by the group ${\rm Aut}\, (X/S)$. However, the class
$K_{X/S}$ is relatively trivial over $S$, if $f$ is a relatively minimal Kodaira-Iitaka
fibration or if $X$ has a numerically trivial canonical class ($X$
a Calabi-Yau variety in a wider sense) and $S = {\rm Spec}\, k$.

The next theorem will be applied in Section 5.

\begin{theorem} \label{theorem:chern}
Let $X$ be an $n$-dimensional projective manifold whose $(n-1)$-th
Chern class $c_{n-1}(X)$ is in the boundary of the cone of
effective $1$-cycles $\overline{\NE}(X)$. Let $G$ be a numerically
virtually solvable subgroup of ${\rm Aut}\,(X)$. Then there are a
finite-index subgroup $G_1$ of $G$ and a nonzero real nef vector
$v$ such that $G_1(\mathbf R_{\ge 0}v) = \mathbf R_{\ge 0}v$ and
$(v.c_{n-1}(X)) = 0$. The same is true for a complex projective
minimal threefold $($see $\S 5)$ whose second Chern class
$c_{2}(X)$ is not in the interior of $\overline{\NE}(X)$.
\end{theorem}

\begin{proof} Since $c_{n-1}(X)$ is stable under ${\rm Aut}\,(X)$,
 the free submodule $M$ of $\NS(X)$ defined by
$$M := \{x \in \NS(X)\, \, \vert\, \, (x, c_{n-1}(X)) = 0\}$$
is also ${\rm Aut}\, (X)$-stable. We set $$C := \overline{{\rm
Amp}}\,(X) \cap M_{\mathbf R}.$$ Since $c_{n-1}(X)$ is in the
boundary of $\overline{\NE}(X)$, it follows that $C \not= \{0\}$.
So, we obtain the first assertion by applying Corollary 2.3 for the
action of $G$ on $(M, C)$. The last statement follows from a
result of Miyaoka [Mi], which says that $c_{2}(X)$ is
well-defined and is always in $\overline{\NE}(X)$ for a minimal
projective complex threefold.
\end{proof}

\medskip
{\bf (3.5)} For several reasons, even ${\rm Bir}^{s}(X/S)$ does
not act on $\NS(X/S)$ in general. We quote the following result of
Kawamata [Ka2, the proof of Lemma 1.1].

\begin{theorem} \label{theorem:index}

Assume that ${\rm char}\, k = 0$ and $X$ admits at most $\mathbf
Q$-factorial rational singularities. Let $\N^{1}(X/S)$ be the group
of numerical equivalence classes of Weil divisors on $X$ over $S$.
Then $\N^{1}(X/S)$ is a free $\mathbf Z$-module of finite rank and
satisfies $[\N^{1}(X/S) : \NS(X/S)] < \infty$. Moreover, ${\rm
Bir}^{s}(X/S)$ acts naturally and linearly on $\N^{1}(X/S)$ and
preserves $\overline{{\rm Big}}(X/S)\, (\subset
\N^{1}(X/S)_{\mathbf R})$.
\end{theorem}

Let $G$ be a subgroup of ${\rm Bir}^{s}\,(X/S)$. We call $G$ {\it
numerically virtually solvable} if its image in ${\rm GL}\,
(\N^{1}(X/S))$ is virtually solvable. Since $[\N^{1}(X/S) : \NS(X/S)] <
\infty$, the definitions here and in (3.4) are consistent if $G
\subset {\rm Aut}\, (X/S)$.

Now, applying Corollary 2.3 to the action $G$ on $(\N^{1}(X/S), \,
\overline{{\rm Big}}(X/S))$, we obtain the following:

\begin{theorem} \label{theorem:bir}
Assume that ${\rm char}\, k = 0$ and $X$ admits at most $\mathbf Q$-factorial rational
singularities. Let $G$ be a numerically
virtually solvable subgroup of ${\rm Bir}^{s}\,(X/S)$. Then there
are a finite-index subgroup $G_1$ of $G$ and a nonzero real vector
$v \in \overline{{\rm Big}}(X/S)$ such that $G_1(\mathbf R_{\ge
0}v) = \mathbf R_{\ge 0}v$. In particular, if $X$ is minimal and
if $G$ is a numerically virtually solvable subgroup of ${\rm
Bir}\, (X)$, then a suitable finite-index subgroup of $G$ admits a
common real pseudo-effective eigenvector.
\end{theorem}

In the above Theorem, we do not need the rationality of
singularities, because $\mathbf Q$-factoriality already implies
that $\NS (X/S)_{\mathbf Q} = \N^{1}(X/S)_{\mathbf Q}$. However,
in order to guarantee the discreteness as in Lemma 2.4, the action
of $G$ on some $\mathbf Z$-module structure will be needed.

\medskip
 {\bf (3.6)} Let $X$ be a compact K\"ahler manifold and let
$\overline{\mathcal K}(X)$ be the closure of the K\"ahler cone
$\mathcal K(X)$ in $H^{2}(X, \mathbf R)$. Let $G$ be a subgroup of
${\rm Aut}\, (X)$. We call $G$ {\it cohomologically virtually
solvable} if its image in ${\rm GL}\, (H^{2}(X, \mathbf Z)/ {\rm
torsion})$ is virtually solvable. When $X$ is projective, $G$ is
numerically virtually solvable if $G$ is cohomologically virtually
solvable.

Applying Corollary 2.3 to the action of $G$ on $(H^{2}(X, \mathbf
Z)/ {\rm torsion}, \, \overline{\mathcal K}(X))$, we obtain the
following:

\begin{theorem} \label{theorem:kahler}
If $G$ is a cohomologically virtually solvable subgroup of ${\rm
Aut}\, (X)$, then there are a finite-index subgroup $G_1$ of $G$
and a nonzero real vector $v \in \overline{\mathcal K}(X)$ such
that $G_1(\mathbf R_{\ge 0}v) = \mathbf R_{\ge 0}v$.
\end{theorem}

\section{Entropy and Conjecture of Tits type}

In this section, we verify Conjecture of Tits type (Conjecture
1.3) posed in Introduction for some basic varieties, i.e., for
surfaces, complex tori and hyperk\"ahler manifolds. In the next
section, we shall also verify the conjecture for minimal
threefolds. In our proof, Theorem 2.1 or Corollary 2.3
also plays an important
role.

\medskip
{\bf (4.1)} We recall a definition of {\it entropy}
 of an automorphism of a compact K\"ahler manifold from [DS] and
generalize it to a complex normal projective variety.

For the references to the following important results,
see Dinh-Sibony [DS] and Guedj\cite[(1.2), (1.6)]{Gu}.

\begin{theorem} \label{theorem:ent} $([DS])$
Let $M$ be a compact K\"ahler manifold of dimension $n$ and let
$g$ be an automorphism of $M$. By $\rho(g^{*} \vert W)$ we denote
the spectral radius of the action of $g^*$ on a $g^*$-stable
subspace $W$ of the total cohomology group $H^{*}(M, \mathbf C)$.
Then we have:

\begin{itemize}
\item[(1)] $\rho(g^{*} \vert H^{*}(M, \mathbf C)) \ge 1$, and $\rho(g^{*}
\vert H^{*}(M, \mathbf C)) = 1$ (resp. $> 1$) if and only if
$\rho(g^{*} \vert H^{2}(M, \mathbf C)) = 1$ (resp. $> 1)$.

Moreover, $\rho(g^{*} \vert H^{*}(M, \mathbf C)) = 1$
(resp. $> 1$) if and only if so is for $g^{-1}$.

\item[(2)] $\rho(g^{*} \vert H^{2}(M, \mathbf C)) = \rho(g^{*} \vert H^{1,1}(M))$.

Moreover, if $M$ is projective, then this value is also equal to
$\rho(g^{*} \vert \NS(M))$.

\item[(3)] The spectral radius $\rho(g^*|H^{i,i}(X, {\mathbf C}))$ equals
the $i$-th dynamical degree $d_i(g)$. Further, there are integers
$m \le m'$ such that
$$1 = d_0(g) < \dots < d_m(g) = \dots = d_{m'}(g)
> \dots > d_n(g) = 1.$$

\item[(4)]
$\rho(g^* | H^*(X, {\mathbf C})) = \max_{0 \le i \le n} \, d_i(g)
= e^{h(g)}$, with $h(g)$ the topological entropy of $g$.

\end{itemize}
\end{theorem}

\begin{definition} \label{definition:entropy}
$ $
\par
(1) Let $M$ be a compact K\"ahler manifold and let $g$ be an
automorphism of $M$. We call $g$ {\it of null entropy} (resp. {\it
of positive entropy}) if and only if the spectral radius
$\rho(g^{*} \vert H^{2}(M, \mathbf C)) = 1$ (resp. $> 1$). By
(4.1), $g$ is of null entropy (resp. of positive entropy) if and
only if $\rho(g^{*} \vert \NS(M)) = 1$ (resp. $> 1$) when $X$ is
projective.

(2) Let $V$ be a complex normal projective variety. Taking account
of the last characterization in (1), we call an automorphism $g$
of $V$ {\it of null entropy} (resp. {\it of positive entropy}) if
the spectral radius $\rho(g^{*} \vert \NS(V))$ on the free part
$\NS(V)$ of the N\'eron-Severi group is $1$ (resp. $> 1$).
If $V$ is {\bf Q}-factorial, then $g \in {\rm Aut}(V)$ is of positive entropy
if and only if so is the lifting of $g$ to an automorphism
on an equivariant resolution of $V$; see the proof of \cite[Lemma 2.6]{Zh2}.
\end{definition}

\medskip
{\bf (4.2)} In Introduction, inspired by Tits' Theorem and
Dinh-Sibony's Theorem (Theorem 1.2) we posed  Conjecture of Tits
type (Conjecture 1.3). The following is a weaker version of the
conjecture, which can be easily verified.

\begin{theorem} \label{theorem:weak}
Let $L$ be a free $\mathbf Z$-module of rank $r \ge 1$ and let $H$ be a
subgroup of ${\rm GL}(L)$. Assume that $H$ does not contain a
subgroup isomorphic to $\mathbf Z * \mathbf Z$. Then, there is a
connected solvable finite-index subgroup $H_1$ of $H$ such that the
subset
$$N(H_1) := \{h \in H_1\, \, \vert\, \, \rho(h) = 1\}$$
of $H_1$ is a normal subgroup of $H_1$ and $H_1/N(H_1)$ is a free
abelian group of rank at most $r-1$.
\end{theorem}

Applying this theorem to the action of $G$ on $H^{2}(X, \mathbf
Z)$ or on $\NS(X)$, we see that the conjecture is true {\it except
the rank estimate}. In other words, {\it the rank estimate is the
most essential part of the conjecture}.
\begin{proof} By Tits' Theorem, $H$ is virtually solvable. Thus $H$ has a
finite-index subgroup $H_1$ which is connected solvable. Then, the
result follows from Lemma 2.4.
\end{proof}

\medskip
{\bf (4.3)} In the rest of this section, we verify Conjecture of
Tits' type for surfaces (Theorem 4.4), hyperk\"ahler manifolds
(Theorem 4.5), and complex tori (Theorem 4.6). We also study
Conjecture of Tits type for coverings and fibrations (Propositions
4.7, 4.8).

\begin{theorem} \label{theorem:surface}
Conjecture of Tits type is true if $\dim X \le 2$.
\end{theorem}

\begin{proof} There is nothing to prove
when $\dim X = 0$ or $1$.

Assume $\dim X = 2$. If $X$ is smooth and projective, the
conjecture has been proved by [Zh1]. Therefore, when $X$ is
projective, it is possible to reduce to the smooth projective case
by considering a minimal resolution $\tilde{X}$ of $X$ and the
induced action on $\tilde{X}$. It is also possible to argue by
using classification of surfaces (see eg. [BHPV]) when $X$ is
K\"ahler. Instead of such a case-by-case proof, we shall give here
a unified, classification-free proof.

Define $$L:=H^{2}(X, \mathbf Z)\,\,{\rm (resp.} =\NS(X))\,,$$
$$C:=\overline{\mathcal K}(X)\,\,{\rm (resp.} =\overline{{\rm
Amp}}(X))$$ when $X$ is a compact K\"ahler surface (resp. $X$ is a
normal projective surface).
%
%
%
%
By the Hodge index theorem, the real
linear subspace $$W := \langle C \rangle$$ of $L_{\mathbf R}$ is
hyperbolic, i.e. of signature $(1, *)$ with respect to the
intersection form. This is because $$W = H^{1,1}(M, \mathbf
R)\,\,{\rm or}\,\, W = \NS(X)_{\mathbf R}$$ according to the two
cases above.

Let us consider the image
$$G^{*} := {\rm Im}\, (G \longrightarrow {\rm GL}(L))\,\, .$$
We may assume that $G^{*}$ is connected solvable. Let $v \in C$ be
a common eigenvector of $G^{*}$ found by Corollary 2.3. Then, we can
write $g^{*}(v) = \chi(g)v$.

\medskip
{\bf Claim :} (1) $\rho(g^*\vert W) = \rho(g^{*} \vert L)$.
\,\,\,\,\,\,(2) $\chi(g) = \rho(g^{*} \vert L)\,\,\,\, {\rm or}\,\,\,\,
\frac{1}{\rho(g^{*} \vert L)}\,\,.$

\medskip
{\it Proof of Claim.}  (1) In the K\"ahler surface case, the
intersection form on the subspace $(H^{2,0}(X)\oplus
H^{0,2}(X))\cap H^{2}(X, {\mathbf R})$ is positive definite, so
the eigenvalues of $g^*$ on this subspace all have absolute value
$1$.

(2)  By Birkhoff-Perron-Frobenius' Theorem (Remark 2.2(3)) applied
to the action of $g^*$ on $(W, C)$, there is a nonzero vector
$v_{g}\in C$ such that $g^{*}(v_{g}) = \rho(g^*\vert W) \ v_{g}$.
Note that
$$(v. v_{g}) = (g^{*}(v). g^{*}(v_{g}))
= \chi(g)\rho(g^*\vert W)(v. v_{g})\,\, .$$ If $\chi(g)\neq
\frac{1}{\rho(g^*\vert W)}$\,\,, then we must have $(v, v_{g}) =
0$. Since $v, \,  v_{g} \in C$ and since $W$ is hyperbolic by the
Hodge index theorem, it follows that $v$ is parallel to $v_{g}$,
so $\chi(g)=\rho(g^*\vert W)$. Now, (2) follows from (1).

\medskip
Resuming the proof of the theorem, we consider the group
homomorphism
$$\chi : G \longrightarrow \mathbf R\,\, ;\,\, g \mapsto
\log\vert \chi(g) \vert=\log \chi(g)\,\, .$$ By Claim, we see that
$\chi(g) = \rho(g^{*} \vert L)=1$ if and only if $g$ is of null
entropy. This means that $${\rm Ker}\, \chi = N(G),$$ and
hence
$$G/N(G) \simeq \chi(G).$$
Furthermore, by Claim again, we have
$$\chi(G) \subset \{\pm\log \rho(g^{*} \vert L)\,\, \vert\,\,  g \in G\}\,\, .$$
Note here that if $\rho(g^{*}\vert L)$ is bounded, then all
eigenvalues of $g^{*} \vert L$ are also bounded. We also note that
the action $G^{*}$ is defined over the $\mathbf Z$-module $L$.
Thus, $\chi(G)$ is discrete in $\mathbf R$ as we have seen in the
proof of Lemma 2.4. Hence $\chi(G) \simeq \mathbf Z^{s}$ with $s
\le 1$.
\end{proof}

By a {\it hyperk\"ahler} manifold $M$, we mean a simply-connected
compact K\"ahler manifold admitting an everywhere non-degenerate
holomorphic $2$-form $\sigma_{M}$ such that $H^{0}(M,
\Omega_{M}^{2}) = \mathbf C \sigma_{M}$. According to the
Bogomolov decomposition Theorem [Be1], hyperk\"ahler manifolds
form one of the three building blocks of compact K\"ahler
manifolds with vanishing first Chern class. The other two building
blocks are complex tori and Calabi-Yau manifolds (in the narrow
sense).

\begin{theorem} \label{theorem:hk}
Conjecture of Tits type is true for hyperk\"ahler manifolds.
\end{theorem}

\begin{proof} By Beauville [Be1], $H^{2}(M, \mathbf Z)$ of a
hyperk\"ahler manifold $M$ admits an integral symmetric bilinear
form $(*, **)$ of signature $(3, b_{2}(M) -3)$. This form is
invariant under ${\rm Aut}\,(M)$ and is of signature $(1,
h^{1,1}(M)-1)$ on $H^{1,1}(M, \mathbf R)$. A bit more precisely,
$(\eta, \eta') > 0$ for any K\"ahler classes $\eta$ and $\eta'$.
Thus, the same proof as in Theorem 4.4 shows the result.
Furthermore, in the statement (2) of the conjecture, we have
$G_1/N(G_1) \simeq \mathbf Z^{s}$ with $s \le 1$.
\end{proof}

\begin{theorem} \label{theorem:tori}
Conjecture of Tits type is true for complex tori.
\end{theorem}

\begin{proof} Let $X$ be a complex torus of dimension $n$ and
let $G$ be a subgroup of ${\rm Aut}\, (X)$. Note that
$$H^{1}(X, \mathbf Z) \otimes \mathbf C =
H^{0}(X, \Omega_{X}^{1}) \oplus \overline{H^{0}(X, \Omega_{X}^{1})}\,\, ,$$
$$H^{*}(X, \mathbf Z) = \oplus_{k=0}^{2n} \wedge^{k} H^{1}(X, \mathbf Z)\,\, .$$
Thus $g \in G$ is of null entropy if and only if the spectral radius
$\rho(g^{*} \vert H^{0}(X, \Omega_{X}^{1})) = 1$. Let us consider
the action of $G$ on $H^{0}(X, \Omega_{X}^{1}) \simeq \mathbf
C^{n}$, and define
$$G^{*} := {\rm Im}\, (G \longrightarrow {\rm GL}\, (H^{0}(X, \Omega_{X}^{1}))\,\, .$$
As in the proof of Theorem 4.4, we may assume that $G^{*}$ is
connected solvable. Then, as in Lemma 2.4, choosing a suitable
basis of $H^{0}(X, \Omega_{X}^{1})$, one can embed $G^{*}$ into
the subgroup $T(n)$ of upper triangle matrices of ${\rm GL}\, (n,
\mathbf C) = {\rm GL}\,(H^{0}(X, \Omega_{X}^{1}))$, and we have a
group homomorphism
$$\chi : G \longrightarrow \mathbf R^{n-1} = \{(x_{i})_{i=1}^{n}
\in \mathbf R^{n} \, \vert \, \sum_{i=1}^{n} x_{i} = 0\} \subset \mathbf R^{n}$$
which is defined by
$$\chi(g) = (\log \vert \chi_{1}(g) \vert, \ \log \vert \chi_{2}(g)\vert, \
\cdots , \ \log \vert \chi_{n}(g)\vert)\,\, .$$
Here $\chi_{i}(g)$ is the $(i,i)$-th entry of the matrix $g$
(regarded as an element of $T(n)$). Thus $$N(G) = {\rm Ker}\,
\chi\,.$$ Note that the eigenvalues of $g^{*} \vert H^{1}(X,
\mathbf Z)$ are $\chi_{i}(g)$ and $\overline{\chi_{i}(g)}$ ($1 \le
i \le n$). Thus, one can apply the same argument as in Lemma 2.4
to get $$G/N(G) \simeq {\rm Im}\,\chi \simeq \mathbf Z^{s}$$ for
some $s \le n-1$.
\end{proof}

It is highly interesting to verify the conjecture for Calabi-Yau
manifolds (in the narrow sense). Unfortunately, we could not yet
do it. All we can say now is that the conjecture is true for
Calabi-Yau threefold, as a special case of Theorem 5.1 in Section
5. See also Proposition 4.9 for a relevant result.

Next two propositions will be used in Section 5.

\begin{proposition} \label{proposition:covering}
Let $G$ be a subgroup of ${\rm Aut}\, X$ and let $\pi : \tilde{X}
\longrightarrow X$ be a generically finite surjective morphism
between compact K\"ahler manifolds or between complex normal
projective varieties. Assume that $G$ lifts to a subgroup
$\tilde{G} \simeq G$ of ${\rm Aut}\,(\tilde{X})$ equivariantly.
Then Conjecture of Tits type is true for $(X, G)$ if so is for
$(\tilde{X}, \tilde{G})$.
\end{proposition}

\begin{proof} Let $\dim X = \dim \tilde{X} = n$. If $G$ contains a
subgroup isomorphic to $\mathbf Z * \mathbf Z$, then we are done.
Therefore we may assume that $\tilde{G} \simeq G$ does not contain
a subgroup isomorphic to $\mathbf Z * \mathbf Z$. Then, by Tits'
Theorem, any linear action of $\tilde{G}$ on a free $\mathbf
Z$-module of finite rank is virtually solvable. Furthermore, by
replacing $G \simeq \tilde{G}$ by a suitable finite-index
subgroup, we may assume that $\tilde{G} \vert H^{2}(\tilde{X},
\mathbf Z)$ and $G \vert H^{2}(X, \mathbf Z)$ or $\tilde{G} \vert
\NS(\tilde{X})$ and $G \vert \NS(X)$ are connected solvable. Since
$\pi^{*} : H^{2}(X, \mathbf Z) \longrightarrow H^{2}(\tilde{X},
\mathbf Z)$ or $\pi^{*} : \NS(X) \longrightarrow \NS(\tilde{X})$
is injective, there is then a surjective homomorphism
$\phi:\tilde{G}/N(\tilde{G}) \longrightarrow G/N(G)$. Since
$G/N(G)$ is free by Lemma 2.4, the result follows.
\end{proof}

If both $X$ and $\tilde{X}$ are {\bf Q}-factorial, then the group
homomorphism $\phi$ is indeed an isomorphism ( see \cite[Appendix,
Lemma A.8]{NZ} and also \cite[Lemma 2.6]{Zh2}), hence the converse
of the conclusion is also true.

Thanks to Proposition 4.7, we may freely replace the variety $X$
by its $G$-equivariant resolution, which always exists.

\begin{proposition} \label{proposition:fibration}
Let $f : X \longrightarrow Y$ be a surjective morphism with
connected fibres. Assume that $X$ is smooth. Let $F$ be a general
fiber. Let $G$ be a subgroup of ${\rm Aut}\, (X/Y)$ $($so that $G$
acts faithfully on $F )$. Then, Conjecture of Tits type is true
for $G$ with the rank bound $\dim F -1$ if and only if so is for
the action of $G$ on $F$.
\end{proposition}

\begin{proof} We write by $G_{F}$ the action of $G$ on $F$.
As in the previous proposition, we may assume that $G$ does not
contain ${\mathbf Z} * {\mathbf Z}$. Replacing $G$ by a
finite-index subgroup, we may assume that both
$$G^{*} := {\rm Im}(\varphi_{X} : G \longrightarrow {\rm
GL}\,(\NS(X))\,\,{\rm and}\,\,G_{F}^{*} := {\rm Im}(\varphi_{F} :
G_{F} \longrightarrow {\rm GL}\,(\NS(F))$$ are connected solvable.
The groups $N(G)$ and $N(G_{F})$ are then well-defined normal
subgroups of $G = G_{F}$ (cf. Lemma 2.4). Then, by [Zh2] (2.1)
Remark (11) or \cite[Appendix, Theorem D]{NZ}, we have $N(G_{F}) =
N(G)$ as subgroups of $G$. Thus we have the identification:
$$G_{F}/N(G_{F}) =  G/N(G)\,\,.$$
This implies the result.
\end{proof}

As an application, we have the following:

\begin{proposition} \label{proposition:trivialchern}
Conjecture of Tits type is true for compact K\"ahler manifolds $M$
with $c_{1}(M) = 0$ in $H^{2}(M, \mathbf R)$, provided that
Conjecture of Tits type is true for Calabi-Yau manifolds $($in the
narrow sense$)$.
\end{proposition}

\begin{proof} Let $\pi : \tilde{M} \longrightarrow M$
be the minimal splitting cover of $M$ in the sense of [Be2].
$\tilde{M}$ is the product of complex torus $T$, finitely many
hyperk\"ahler manifolds $V_{i}$ and finitely many Calabi-Yau
manifolds $W_{j}$. The subgroup $G$ of ${\rm Aut}\, (M)$ lifts
equivariantly to the subgroup $\tilde{G} \simeq G$ of ${\rm Aut}\,
(\tilde{M})$ and $\tilde{G}$ preserves each factor $T$, $V_{i}$,
$W_{j}$ [ibid]. By the K\"unneth formula, $H^{2}(\tilde{M},
\mathbf Z)$ is naturally isomorphic to the direct sum of $H^{2}(T,
\mathbf Z)$, $H^{2}(V_{i}, \mathbf Z)$ and $H^{2}(W_{j}, \mathbf
Z)$. Thus, Conjecture of Tits type is true for $\tilde{G}$ by
Theorems 4.5, 4.6 and by our assumption for Calabi-Yau manifolds.
Now the result follows from Proposition 4.7.
\end{proof}

\begin{remark} \label{remark:noncy}
Furthermore, the proof of Proposition 4.9 together with Theorem
5.1 below implies that Conjecture of Tits type is true for compact
K\"ahler manifolds $M$ with $c_{1}(M) = 0$ in $H^{2}(M, \mathbf
R)$ if its minimal splitting cover does not have a Calabi-Yau
manifold of dimension $\ge 4$ as a factor.
\end{remark}

\section{Conjecture of Tits type for minimal threefolds}

Let $X$ be a normal complex projective variety. It is called
{\it minimal} if
$X$ admits at most $\mathbf Q$-factorial terminal singularities
and the canonical class $K_X$ is nef.

In this section, we shall show the following:

\begin{theorem} \label{theorem:threefold}
Conjecture of Tits type is true for minimal threefolds.
\end{theorem}

Our proof is mostly classification free, but unfortunately, we
could not avoid using some classification result of Calabi-Yau
threefolds of very special type at the final step.
\begin{proof} By the abundance theorem for threefolds [Ka3],
the linear system $\vert mK_X \vert$ for some $m > 0$ defines a
surjective morphism $\varphi : X \longrightarrow W$ from $X$ to
the canonical model $W$ of $X$ (Iitaka-Kodaira fibration). The
morphism $\varphi$ is ${\rm Aut}\,(X)$-equivariant and the induced
action of ${\rm Aut}\, X$ on $W$ is finite by
Deligne-Nakamura-Ueno \cite[Theorem 14.10]{Ue}. Thus, replacing
$G$ by a finite-index subgroup, we may assume that $G \subset {\rm
Aut}\, (X/W)$. Since the conjecture of Tits type is true when
$\dim \le 2$ by Theorem 4.4, the conjecture is also true if $\dim
W \ge 1$, i.e. if $K_{X}$ is not numerically trivial. This can be
seen by applying Proposition 4.7 to a $G$-equivariant resolution
$\pi:\tilde{X}\longrightarrow X$, and then Proposition 4.8 to the
composition $\varphi\circ\pi : \tilde{X} \longrightarrow W$.

Let us consider the case where $K_{X}$ is numerically trivial. By
Proposition 4.7 applied to the global index-one cover of $X$, we
may assume that $K_{X}$ is linearly equivalent to $0$.

If $q(X) := h^{1}(\mathcal O_{X}) > 0$, then the albanese morphism
$a : X \longrightarrow \Alb(X)$ is an \'etale fiber bundle by
[Ka1]. Thus, $X$ is smooth and its minimal splitting cover is
either the product of an elliptic curve and a K3 surface or an
abelian threefold. In this case, the result follows from
Proposition \ref{proposition:trivialchern}.

If $c_{2}(X) = 0$ as a linear form on $\NS(X)$, then $X$ is an
\'etale quotient of an abelian threefold by [SW]. So, $X$ is
smooth and its minimal splitting cover is an abelian threefold.
Therefore the result also follows from Proposition
\ref{proposition:trivialchern}.

Thus, we may assume further that $q(X) = 0$ and $c_{2}(X) \not= 0$.

\begin{lemma} \label{lemma:either}
Let $X$ be a minimal threefold such that $K_X \sim 0$, $q(X) = 0$
and $c_{2}(X) \not= 0$ $($or more generally $K_X \sim_{\mathbf Q}
0$ and $c_{2}(X) \not= 0)$. Let $G \subset {\rm Aut}\, (X)$. Then,
either one of the following holds:

\begin{itemize}
\item[(1)] Conjecture of Tits type is true for $(X, G)$.
\item[(2)] There is a nef and big Cartier divisor $D$ such that
$(c_{2}(X).D) = 0$.
\end{itemize}
\end{lemma}

\begin{proof} Set $G^{*} := {\rm Im}\, (G \longrightarrow {\rm GL}\, (\NS(X)))$.
If $\mathbf Z * \mathbf Z \subset G^{*}$, then (1) holds.
So, by Tits' theorem, we may assume that $G^{*}$ is connected solvable
(after replacing $G$ by a finite-index subgroup).

Let
$$M := \{x \in \NS(X)\,\, \vert\,\, (c_{2}(X).x) = 0\,\}\,\, ,\,\,\,\,
C := \overline{{\rm Amp}}\, (X) \cap M_{\mathbf R}\,\, .$$ For
$g\in G$, we denote by $\rho(g)$ the spectral radius of $g^{*}
\vert \NS(X)$, i.e. $$\rho(g):= \rho(g^{*} \vert \NS(X))\,.$$ If
each element of $G$ is of null entropy, then we are done. So we
may assume that $\rho(g_1) > 1$ for some $g_1 \in G$. Then, by
Birkhoff-Perron-Frobenius' Theorem, there is a real vector
$v_{g_1} \in \overline{{\rm Amp}}\, (X) \setminus \{0\}$ such that
$g_1^*v_{g_1}= \rho(g_1)v_{g_1}$. Thus $c_2(X) . v_{g_1} = 0$ as
shown below. Hence $C \neq \{0\}$ and we can apply Theorem 3.2 to
obtain a real vector $v_{0} \in C \setminus \{0\}$ such that
$G(\mathbf R_{\ge 0}v_{0}) = \mathbf R_{\ge 0}v_{0}$. Thus, for
each $g\in G$ we can write $g^{*}v_{0} = \chi(g)v_{0}$, and hence
we have the group homomorphism:
$$\chi : G \longrightarrow \mathbf R\,\, ;\,\, g \mapsto \log \chi(g)\,\, .$$

 Note
that $\rho(g) > 1$ if and only if $\rho(g^{-1})
> 1$ (for ${\rm det} \,  g^{*} \vert \NS(X) = \pm 1$). Furthermore, if $\rho(g)=1$ then
$\rho(g)=\chi(g)=\rho(g^{-1})=1$. This implies that we have either
one of the following two cases:

\begin{itemize}
\item[Case(i)] There is an element $g\in G$ such that $\rho(g)
> 1$ and the three positive real numbers
$\chi(g),\,\, \rho(g),\,\,1/\rho(g^{-1})$ are mutually distinct.
\item[Case(ii)] For each $g\in G$, either $\rho(g)=1$, or
$\chi(g)=\rho(g)>1$, or $\chi(g)=1/\rho(g^{-1})<1$.
\end{itemize}
\medskip

We deal with Case (i) first.
By Birkhoff-Perron-Frobenius' Theorem, there is a
real vector $v_{g} \in \overline{{\rm Amp}}\, (X) \setminus \{0\}$
such that
$$g^{*}v_{g} = \rho(g) \ v_{g}\,.$$ Since $g^{*}c_{2}(X) = c_{2}(X)$,
we have
$$(c_{2}(X).v_{g}) = (g^{*}c_{2}(X).g^{*}v_{g}) = \rho(g)(c_{2}(X).v_{g})\,\, .$$
Thus $$(c_{2}(X).v_{g}) = 0\,.$$ Similarly,
since $\rho(g^{-1}) > 1$, there is a real vector $v_{g^{-1}} \in
\overline{{\rm Amp}}\, (X) \setminus \{0\}$ such that
$$(g^{-1})^*v_{g^{-1}} = \rho(g^{-1})v_{g^{-1}}\,\,\,\,{\rm and}
\,\,\,\,(c_{2}(X).v_{g^{-1}}) = 0\,.$$
Note that $c_2(X)$ gives an integer valued linear form on
$\NS(X)$, so $v_{g}, v_{g^{-1}} \in M_{\mathbf R}$.
Hence
$$v_{g}, \, v_{g^{-1}} \in C \setminus \{0\}\,.$$
Let
$$D' := v_{0} + v_{g}+v_{g^{-1}}\,.$$ Then $$(D'. c_{2}(X)) = 0\,.$$
Since $v_{0}$, $v_{g}$ and $v_{g^{-1}}$ are eigenvectors of $g^*$
with distinct eigenvalues $\chi(g)$, $\rho(g)$ and
$1/\rho(g^{-1})$, we have $$(v_{0}.v_{g}.v_{g^{-1}}) \not= 0$$ by
[DS], Lemma 4.4. Since $v_{0}$, $v_{g}$ and $v_{g^{-1}}$ are all
nef, we have then $$((D')^{3}) > 0$$ and hence $D'$ is nef and
big. Since $M_{\mathbf R}$ is a rational hypersurface in
$\NS(X)_{\mathbf R}$ and the nef cone is locally rational
polyhedral away from the cubic cone by [Ka2] or [Wi], we can then
find a rational nef and big divisor $D$ such that $(D. c_{2}(X)) =
0$. Therefore, the assertion (2) in the lemma holds.

\medskip

Next we deal with Case (ii).
First note that the image ${\rm Im}\, \chi$ of the
group homomorphism $\chi$ defined above is discrete in $\mathbf
R$. This follows from the same argument as in the last part of the
proof of Theorem 4.4. Thus $${\rm Im}\,\chi \simeq \mathbf
Z^{s}\,\,{\rm with}\,\, s \le 1\,.$$ Also, by the case assumption
it is easy to see that ${\rm Ker}\chi= N(G)$. Therefore,
$G/N(G)\simeq {\rm Im}\,\chi$ and the assertion (1) in the lemma holds.
This proves the lemma.
\end{proof}

Now the following lemma completes the proof of Theorem 5.1.

\begin{lemma} \label{lemma:final}
Let $X$ be a minimal threefold such that $K_{X} \sim 0$,
$q(X) = 0$ and $c_{2}(X) \not= 0$. Assume that there is a nef and big
Cartier divisor $D$ such that $(c_{2}(X).D) = 0$. Then, Conjecture
of Tits type is true for $X$.
\end{lemma}

\begin{proof} Let us write $\zeta_{n} := e^{2\pi i/n}$. Let $l = 3$ or $7$.
By [Og1] or [OS], the universal cover of $X$ is either one of
$X_{l}$ ($l = 3$, $7$). Here $X_{3}$ is the unique crepant
resolution of the quotient threefold $\overline{X}_{3} :=
A_{3}/\langle g_{3} \rangle$, where $A_{3}$ is the $3$-times self
product of the elliptic curve $E_{3}$ of period $\zeta_{3}$ and
$g_{3} := {\rm diag}\, (\zeta_{3}, \zeta_{3}, \zeta_{3})$; $X_{7}$
is the unique crepant resolution of the quotient threefold
$\overline{X}_{7} := A_{7}/\langle g_{7} \rangle$, where $A_{7}$
is the Jacobian threefold of the Klein quartic curve and $g_{7} :=
{\rm diag}\, (\zeta_{7}, \zeta_{7}^{2}, \zeta_{7}^{4})$. Note that
the singular locus of $\overline{X}_{l}$ consists of finitely many
points.

By Proposition 4.7, we may assume that $X$ is $X_{l}$. In each
case, we denote by $\nu_{l} : X_{l} \longrightarrow
\overline{X}_{l}$ the unique crepant resolution. By [OS], this
$\nu_{l}$ is also the unique birational contraction which
contracts the class $c_{2}(X_{l})$ in the cone of effective
$1$-cycles $\overline{\NE}(X_l)$. Thus ${\rm Aut}\, (X_{l})$ acts
equivariantly on $\overline{X}_{l}$

Let $G \subset {\rm Aut}\, (X_{l})$. Then $G$ acts on
$\overline{X}_{l}$ equivariantly. Replacing $G$ by a finite-index
subgroup, we may assume that each $\nu_{l}$-exceptional divisor is
$G$-stable. Thus, it suffices to show the assertion for the action
of $G$ on $\overline{X}_{l}$. Let $\overline{X}_{l}^{0}$ be the
smooth locus of $\overline{X}_{l}$. Clearly $\overline{X}_{l}^{0}$
is $G$-stable. Moreover, by the shape of $g_{l}$, we see that
$\pi_{1}(\overline{X}_{l}^{0}) = \pi_{1}(A_{l}).\langle g_{l}
\rangle$ (a semi-direct product). We regard
$\pi_{1}(\overline{X}_{l}^{0})$ as a deck transformation group of
the universal cover $U_{l} = \mathbf C^{3} \setminus B_{l}$ of
$\overline{X}_{l}^{0}$. Here $B_{l}$ is a discrete set of points
('lying over' the singular points of $\overline{X}_{l}$). Then
$\pi_{1}(\overline{X}_{l}^{0})$ is an affine transformation
subgroup of $\mathbf C^{3}$, in which the factor $\pi_{1}(A_{l})$
forms the group of parallel translations and the factor $\langle
g_{l} \rangle$ forms the group of linear transformations defined
by the matrix $g_{l}$. From this, we can see that, in
$\pi_{1}(\overline{X}_{l}^{0})$,
$$\pi_{1}(A_{l}) = \{\sigma \in \pi_{1}(\overline{X}_{l}^{0})\,\, \vert\,\,
{\rm ord}\, \sigma = \infty\,\}\cup \{\id\}\,\, .$$ Thus the
(equivariant) action of $G$ on the universal cover $U_{l}$
normalizes not only $\pi_{1}(\overline{X}_{l}^{0})$, but also
$\pi_{1}(A_{l})$. Hence the action of $G$ on $\mathbf C^{3}$
descends to an action on $A_{l}$ equivariantly. Now the result
follows from Theorem 4.6 and Proposition 4.7. \end{proof}

This completes the proof of the theorem.
\end{proof}


\end{document}